\documentclass[sn-mathphys,Numbered]{sn-jnl}

\usepackage{subcaption}
\usepackage{graphicx}%
\usepackage{multirow}%
\usepackage{amsmath,amssymb,amsfonts}%
\usepackage{amsthm}%
\usepackage{mathrsfs}%
\usepackage[title]{appendix}%
\usepackage{xcolor}%
\usepackage{textcomp}%
\usepackage{manyfoot}%
\usepackage{booktabs}%
\usepackage{enumitem}
\usepackage{bbm}
\usepackage{ragged2e,booktabs,tabularx,adjustbox}

\newcommand{\scrF}{\mathscr{F}}

\begin{document}

\title[Article Title]{Enhanced Fractional Fourier Transform (FRFT) scheme based on closed Newton-Cotes rules}


\author[]{\fnm{Aubain} \sur{Nzokem}}\email{hilaire77@gmail.com}





\abstract{The paper improves the accuracy of the one-dimensional fractional Fourier transform (FRFT) by leveraging closed Newton-Cotes quadrature rules. Using the weights derived from the Composite Newton-Cotes rules of order QN, we demonstrate that the FRFT of a QN-long weighted sequence can be expressed as two composites of FRFTs. The first composite consists of an FRFT of a Q-long weighted sequence and an FRFT of an N-long sequence. Similarly, the second composite comprises an FRFT of an N-long weighted sequence and an FRFT of a Q-long sequence. Empirical results suggest that the composite FRFTs exhibit the commutative property and maintain consistency both algebraically and numerically. The proposed composite FRFT approach is applied to the inversion of Fourier and Laplace transforms, where it outperforms both the standard non-weighted FRFT and the Newton-Cotes integration method, though the improvement over the latter is less pronounced.}

\keywords{Fractional Fourier Transform (FRFT), Discrete Fourier Transforms (DFT), Newton-Cotes rules, Variance Gamma (VG) distribution, Generalized Tempered Stable (GTS) distribution}



\maketitle
 \section{Introduction}\label{sec1}
 \noindent
Fractional Fourier transform(FRFT) is an important time-frequency analyzing tool, often used for the numerical evaluation of continuous Fourier and Laplace transforms \cite{bailey1994fast,researchmei2010}. FRFT appears in the mathematical literature as early as 1929 \cite{Yang2004} and generalizes the traditional Fourier transform (FT) based on the idea of fractionalizing the eigenvalues of the FT\cite{researchmei2010}. An impetus for studying the fractional Fourier transform is the existence of the fast FRFT algorithm that is significantly more efficient than the conventional fast Fourier transform (FFT) algorithm \cite{bailey1991fractional,garcia1996fractional}. On the other hand, The Newton–Cotes quadrature rules, named after Isaac Newton and Roger Cotes, are the most common numerical integration schemes\cite{chapra2010numerical} based on evaluating the integrand at equally spaced points using the polynomial interpolation. The idea of combining both schemes comes initially from the fact that the fast FRFT scheme is formulated based on a simple step-function approximation to the integral \cite{bailey1994fast}, and the Filon formula \cite{tuck1967simple} was derived on the assumption that the integrand may be approximated stepwise by parabolas. These approximations of the Fourier integrals are called the Filon-Simpson rule, Filon-trapezoidal rule, and, more general, Filon's method\cite{bailey1994fast}. This paper aims to provide a broader development of the approximation evaluation of the fast FRFT from the Newton-cote rules and show that such approximation can be written as the FRFT of the weighted FRFT. The resulting schemes will be applied to analyze the numerical error of two probability density functions. We organize the paper as follows. Section 2 develops the higher-order composite Newton-Cotes quadrature formula. Section 3 presents the fast FRFT algorithm and combines it with Newton-Cote rules. Section 4 provides two illustrative examples.

\section{Composite Newton-Cotes Quadrature Formulas}
\noindent
The Newton-Cotes rules value the integrand $f$ at equally spaced points $x_i$ over the interval $[a,b]$; where $x_i=a+i\frac{b-a}{M}=a+ih$ with $h=\frac{b-a}{M}$; $M=Q{N}$ and $x_{{Qp} + Q}=x_{Q(p+1)}$ where $Q$ is the number of $h$ within the subinterval $[x_{Qp},x_{Qp+ Q}]$ of interval $[a,b]$.
 \subsection{ Composite  Rules} 
\noindent
To have greater accuracy, the idea of the composite rule is to subdivide the interval [a, b] into smaller intervals like $[x_{Qp},x_{Qp+ Q}]$, applying the quadrature formula in each of these smaller intervals and add up the results to obtain more accurate approximations.
 \begin{equation} 
 \label{eqp1}
{\int_{a}^{b}\!f(x)dx}=\sum_{p=0}^{N- 1}{\int_{x_{Qp}}^{x_{Qp + Q}}\!f(x)\,dx}
\end{equation}
We define the Lagrange basis polynomials over the sub-interval $[x_{Qp},x_{Qp+ Q}]$.  
 \begin{equation} 
  \label{eqp2}
  l_{Qp+j}(x)=\prod_{\substack{i\ne j \\  i=0}}^{Q} \frac{x - x_{Qp+i}}{x_{Qp+j} - x_{Qp+i}}  \quad  l_{j}(x_i)=\delta_{ij}= \left\{
  \begin{array}{lr}
    0 & :  i \ne j\\ 
   1& : i=j
  \end{array}
\right.
\end{equation}
\noindent
The Lagrange Interpolating Polynomial  and the integration can be derived 
 \begin{equation} 
 \label{eqp3}
\resizebox{.92\hsize}{!}{$\tilde{f}(x)=\sum_{j=0}^{Q}\!{f(x_{Qp+j})l_{Qp+j}(x)} \quad \int_{x_{Qp}}^{x_{Qp + Q}}\!{\tilde{f}(x)}dx =\sum_{j=0}^{Q}f(x_{Qp+j})\int_{x_{Qp}}^{x_{Qp + Q}}\!l_{Qp+j}(x)\,dx$}
 \end{equation}
\noindent 
The integration of Lagrange basis polynomials
 \begin{equation} 
 \label{eqp4}
\int_{x_{Qp}}^{x_{Qp + Q}}\!{l_{Qp+j}(x)}dx =\frac{b-a}{M} \frac{(-1)^{(Q-j)}}{j!(Q-j)!} \int_{0}^{Q}{\prod_{\substack{i\ne j \\  i=0}}}^{Q}{(y - i)}dy \end{equation}
We have the Lagrange Interpolating integration  over $[x_{Qp},x_{Qp+ Q}]$
   \begin{equation} 
\label{eqp5}
 \int_{x_{Qp}}^{x_{Qp + Q}}\!{\tilde{f}(x)}dx = \frac{b-a}{M}\sum_{j=0}^{Q} W_{j} f(x_{Qp+j}) \quad \quad W_{j} =\frac{(-1)^{(Q-j)}}{j!(Q-j)!} \int_{0}^{Q}\!{\prod_{\substack{i\ne j \\  i=0}}^{Q}{(y - i)}}dy
 \end{equation}

 \textbf{Proposition 1.1 }\label{lem1} \\
 For $Q$ even, $M=Q{N}$ integer,  and  $f \in \mathcal{C}^{Q+2}([a,b])$, there exists $\eta \in ]a,b[ $  such that   \ \  
   \begin{equation} 
\label{eqp10}
\resizebox{.92\hsize}{!}{${\int_{a}^{b}\!f(x)\,dx} = \frac{b-a}{M}\sum_{p=0}^{\frac{M}{Q} - 1} \sum_{j=0}^{Q} W_{j} f(x_{Qp+j}) + {h^{Q+3}}\frac{N f^{(Q+2)}(\eta)}{{(Q+2)}!}\int_{0}^{Q}\!{\int_{0}^{y}\!{\prod_{i=0}^{Q}{(a-i)}\,da} \,dy}$}
 \end{equation}
With $$W_{j} =\frac{(-1)^{(Q-j)}}{j!(Q-j)!} \int_{0}^{Q}{\prod_{\substack{i\ne j \\  i=0}}^{Q}{(y - i)}}dy $$ \\
 \noindent
 For Proposition 1.1  proof, see  \cite{nzokem_2021, aubain2020}\\
 \noindent
The error analysis of the Newton-Cotes formulas of degree $Q$ in \cite{krylov2006approximate,aler2015,aubain2020, nzokem_2021} shows that the level of accuracy is greater when $Q$ is even. In fact, we have $O(h^{Q+3})$ for $Q$ even, against $O(h^{Q+2})$ for $Q$ odd.
 \subsection{Weights Computation}
\begin{figure}[ht]
\vspace{-0.4cm}
    \centering
\hspace{-0.5cm}
  \begin{subfigure}[b]{0.32\linewidth}
    \includegraphics[width=\linewidth]{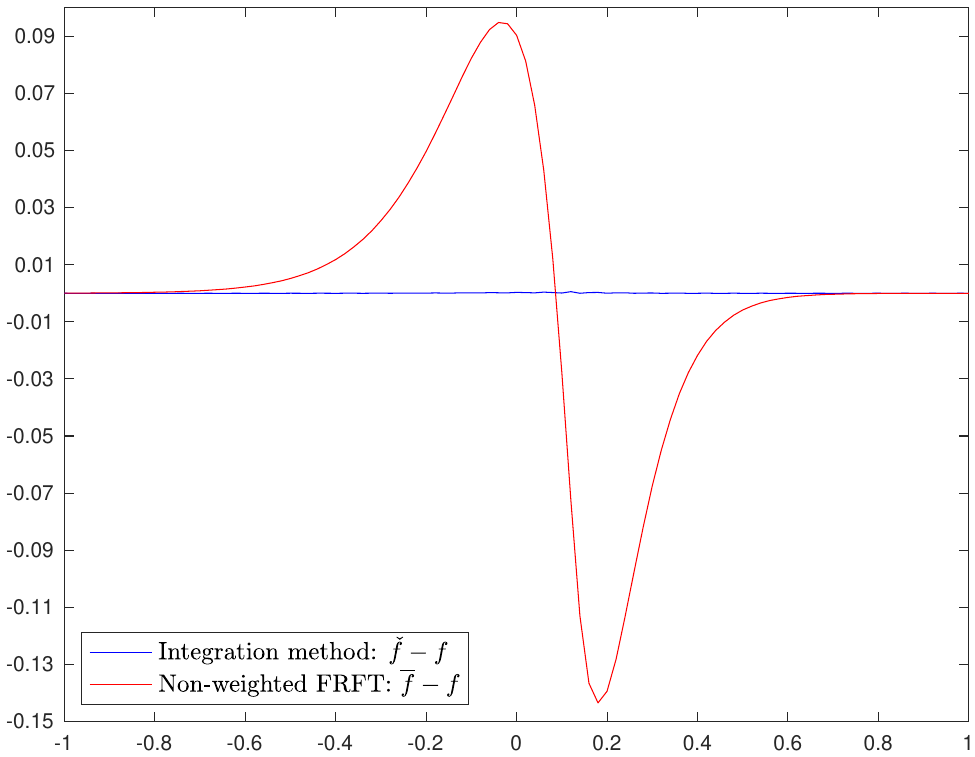}
\vspace{-0.4cm}
     \caption{$\textbf{Q=2}$}
         \label{fig11}
  \end{subfigure}
\hspace{-0.3cm}
  \begin{subfigure}[b]{0.32 \linewidth}
    \includegraphics[width=\linewidth]{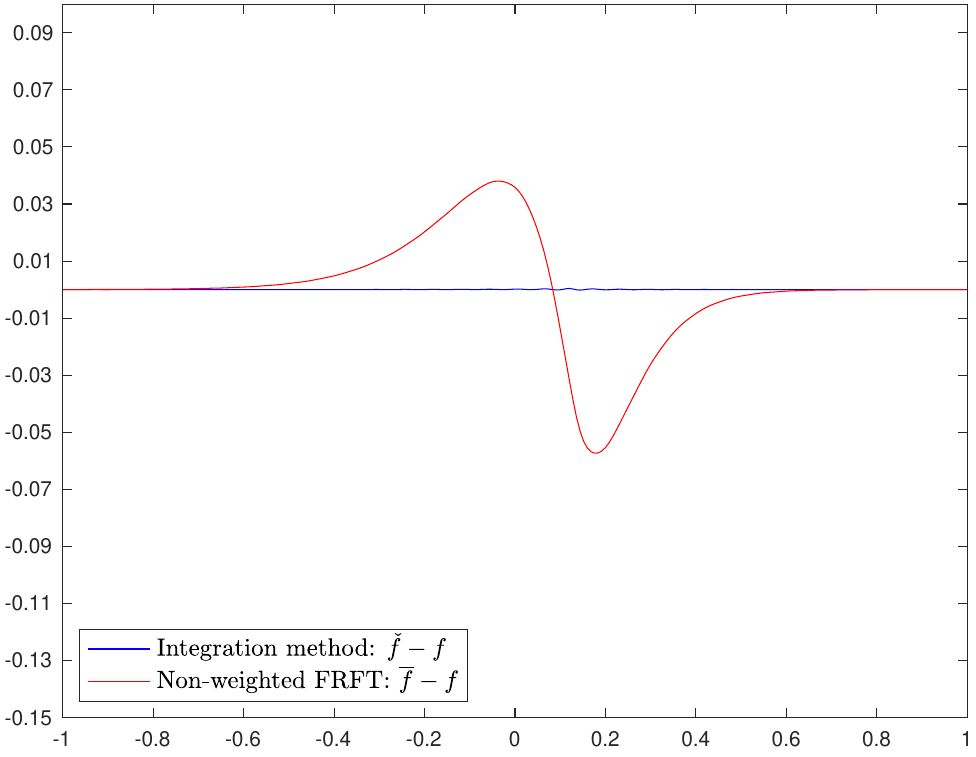}
\vspace{-0.4cm}
     \caption{$\textbf{Q=5}$}
         \label{fig12}
          \end{subfigure}
\hspace{-0.3cm}
  \begin{subfigure}[b]{0.32 \linewidth}
    \includegraphics[width=\linewidth]{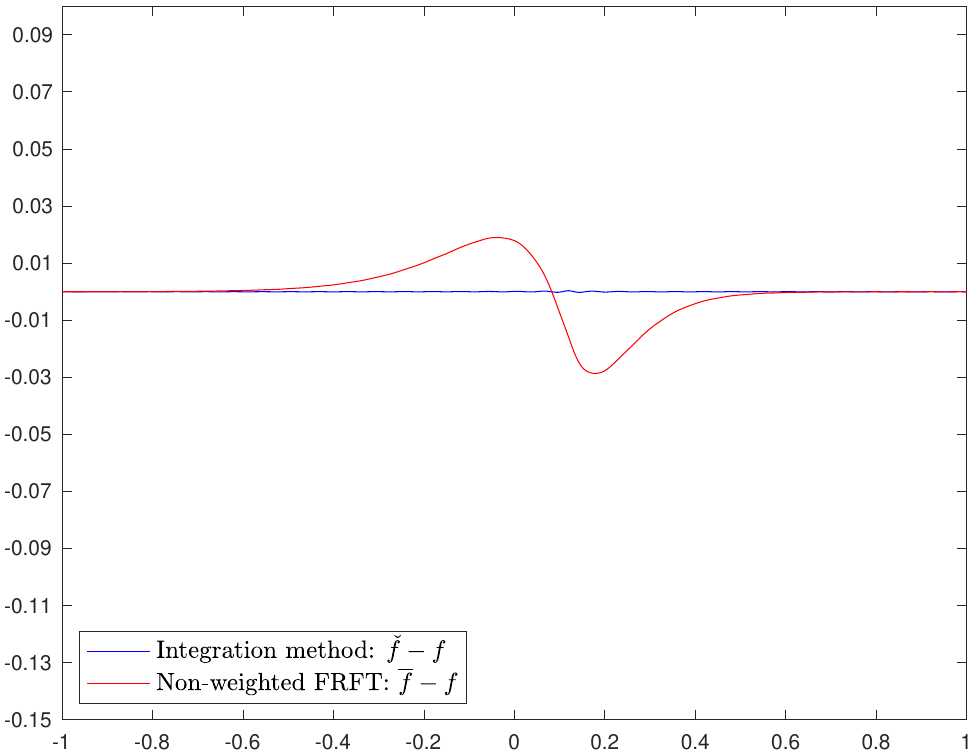}
\vspace{-0.4cm}
     \caption{$\textbf{Q=10}$}
         \label{fig13}
          \end{subfigure}
\vspace{-0.5cm}
  \caption{VG{*} Probability density error ($\textbf{$\check{f}$ - $f$}$): $\textbf{N=5000}$ $\textbf{b-a=100}$}
  \label{fig1}
\vspace{-0.3cm}
\end{figure}

 \noindent
Before using the formula in  (\ref{eqp10}), the expression of  $\{W_j\}_{0\leq j \leq Q}$ can be simplified as follows.\\
\textbf{Proposition 1.2 }\label{lem2} \\
 For $Q$ even and $M=Q{N}$ integer,  $ j\in \{0,1,2 ,...,Q\}$
  \begin{equation} 
  \label{eqp11}
W_{j} = \sum_{i=0}^{Q} C^{j}_{i} \frac{Q^{i+1}}{i+1}\frac{(-1)^{(Q-j)}}{j!(Q-j)!} 
\end{equation} 
 \noindent
 For Proposition 1.2  proof, see  \cite{nzokem_2021,aubain2020}\\
\noindent
The coefficient values $(C^{j}_{i})_{\substack{0\leq i \leq Q  \\ 0\leq j \leq Q}}$ of the polynomial function are obtained by solving the following equations ($\ref{eqp111}$) with a Vandermonde matrix \cite{kalman1984generalized}.
\begin{equation} 
   \label{eqp111}
\prod_{\substack{i\ne j \\ i=0}}^{Q}{(y - i)}=\sum_{i=0}^{Q} C^{j}_{i} y^{i}
\end{equation}
\noindent 
For $Q\leq13$, Table \ref{tab1} provides the values of weights $\{W_j\}_{0\leq j \leq Q}$ in equation (\ref{eqp11}).\\
\noindent
Fig \ref{fig1} compares  the simple non-weighted FRFT and the Composite Newton-Cotes integration method (\ref{eqp10}) using the estimated errors from the risk-neutral probability density of the Variance-Gamma $(\mu, \delta,\alpha,\theta,\sigma)$ model (see  Fig \ref{fig42}). Both algorithms yield different patterns. The Newton-Cotes integration is more accurate.

\begin{table}[ht]
\vspace{-0.3cm}
\centering
\caption{Weights $\{W_j\}_{0\leq j \leq Q}$}
\vspace{-0.3cm}
\label{tab1}
\setlength\extrarowheight{2pt} 
\setlength\tabcolsep{2pt}  
\begin{tabular}{c|c|c|c|c|c|c|c|c|c|c|c|c|c}
\toprule
\textbf{Q} & \textbf{W0} & \textbf{W1} & \textbf{W2} & \textbf{W3} & \textbf{W4} & \textbf{W5} & \textbf{W6} & \textbf{W7} & \textbf{W8} & \textbf{W9} & \textbf{W10} & \textbf{W11} & \textbf{W12}  \\ \midrule
\textbf{1} & $\frac{1}{2}$ & $\frac{1}{2}$ & 0 & 0 & 0 & 0 & 0 & 0 & 0 & 0 & 0 & 0 & 0   \\
\textbf{2} & $\frac{1}{3}$ & $\frac{4}{3}$ & $\frac{1}{3}$ & 0 & 0 & 0 & 0 & 0 & 0 & 0 & 0 & 0 & 0   \\
\textbf{3} & $\frac{3}{8}$ & $\frac{9}{8}$ & $\frac{9}{8}$ & $\frac{3}{8}$ & 0 & 0 & 0 & 0 & 0 & 0 & 0 & 0 & 0   \\
\textbf{4} & $\frac{14}{45}$ & $\frac{64}{45}$ & $\frac{8}{15}$ & $\frac{64}{45}$ & $\frac{14}{45}$ & 0 & 0 & 0 & 0 & 0 & 0 & 0 & 0  \\
\textbf{5} & $\frac{95}{288}$ & $\frac{125}{96}$ & $\frac{125}{144}$ & $\frac{125}{144}$ & $\frac{125}{96}$ & $\frac{95}{288}$ & 0 & 0 & 0 & 0 & 0 & 0 & 0 \\
\textbf{6} & $\frac{41}{140}$ & $\frac{54}{35}$ & $\frac{27}{140}$ & $\frac{68}{35}$ & $\frac{27}{140}$ & $\frac{54}{35}$ &$\frac{41}{140}$ & 0 & 0 & 0 & 0 & 0 & 0 \\
\textbf{7} & $\frac{108}{355}$ & $\frac{810}{559}$ & $\frac{343}{640}$ & $\frac{649}{536}$ & $\frac{649}{536}$ & $\frac{343}{640}$ & $\frac{810}{559}$ & $\frac{108}{355}$ & 0 & 0 & 0 & 0 & 0 \\
\textbf{8} & $\frac{499}{1788}$ & $\frac{1183}{712}$ &$\frac{-182}{695}$ & $\frac{388}{131}$ & $\frac{-319}{249}$ &$\frac{388}{131}$ & $\frac{-182}{695}$ & $\frac{1183}{712}$ & $\frac{499}{1788}$ & 0 & 0 & 0 & 0 \\
\textbf{9} & $\frac{130}{453}$ & $\frac{419}{265}$ & $\frac{23}{212}$ & $\frac{307}{158}$ & $\frac{213}{367}$ & $\frac{213}{367}$ & $\frac{307}{158}$ & $\frac{23}{212}$ & $\frac{419}{265}$ & $\frac{130}{453}$ & 0 & 0 & 0 \\
\textbf{10} & $\frac{139}{518}$ & $\frac{245}{138}$& $\frac{-171}{211}$ & $\frac{414}{91}$ & $\frac{-557}{128}$ & $\frac{1763}{247}$ & $\frac{-557}{128}$ & $\frac{414}{91}$& $\frac{-171}{211}$ & $\frac{245}{138}$& $\frac{139}{518}$ & 0 & 0 \\
\textbf{11} & $\frac{65}{237}$ & $\frac{850}{499}$ & $\frac{-83}{203}$ & $\frac{787}{247}$ &$\frac{-223}{184}$ & $\frac{227}{116}$ & $\frac{227}{116}$ & $\frac{-223}{184}$ & $\frac{787}{247}$ & $\frac{-83}{203}$ & $\frac{850}{499}$ & $\frac{65}{237}$ & 0  \\
\textbf{12} & $\frac{20}{77}$ & $\frac{375}{199}$ & $\frac{-270}{187}$ & $\frac{673}{99}$ & $\frac{-1019}{104}$ & $\frac{816}{49}$ & $\frac{-1537}{92}$ & $\frac{816}{49}$  & $\frac{-1019}{104}$ & $\frac{673}{99}$ & $\frac{-270}{187}$ & $\frac{375}{199}$ & $\frac{20}{77}$  \\ \bottomrule
\end{tabular}
\end{table}

\section{Fast Fractional Fourier Transform (FRFT) and Composite Newton-Cotes Quadrature rules}
\subsection{Fast Fourier Transform and Fractional Fourier Transform }
\noindent
The Conventional fast Fourier transform (FFT) algorithm is widely used to compute discrete convolutions, discrete Fourier transforms (DFT) of sparse sequence, and to perform high-resolution trigonometric interpolation \cite {bailey1991fractional, bailey1994fast}. The discrete Fourier transforms (DFT) are based on $N^{th}$ roots of unity $e^{-\frac{2\pi i}{N}}$. The generalization of DFT is the FRFT, which is based on fractional roots of unity $e^{- 2\pi i\alpha}$, where $\alpha$ is an arbitrary complex number.\\
\noindent
The FRFT is defined on $M$-long sequence ($x_1$, $x_ {2} $, \dots, $x_{M}$) as follows
 \begin{equation}
 G_{k+s}(x,\delta)=\sum_{j=0}^{M-1}\! x_{j}e^{-2\pi i j(k+s)\delta} \hspace{5mm}
   \hbox{$0\leq k<M$ \quad $0\leq s\leq 1$}
   \label {eq:l31}
\end{equation}  
\noindent
Let us have $2j(k+s)=j^2 + (k+s)^2  - (k-j+s)^2$ , equation (\ref{eq:l31}) becomes
 \begin{equation}
 \begin{aligned}
 &G_{k+s}(x,\delta)=\sum_{j=0}^{M-1}\! x_{j}e^{-\pi i (j^2 + (k+s)^2  - (k-j+s)^2) \delta}  \\ \label {eq:l4}
&=e^{-\pi i (k+s)^2 \delta}\sum_{j=0}^{M-1}\! x_{j}e^{-\pi i j^2 \delta}e^{\pi i (k-j+s)^2 \delta}=e^{-\pi i (k+s)^2 }\sum_{j=0}^{M-1}\! y_{j}z_{k-j} 
  \end{aligned}
\end{equation}
\noindent
The expression $ \sum_{j=0}^{M-1}\! y_{j}z_{k-j}$ is a discrete convolution. Still, we need a circular convolution (i.e., $z_{k-j}=z_{k-j+M}$) to evaluate $G_{k+s}(x,\delta)$. The conversion from discrete convolution to discrete circular convolution is possible by extending the sequence $y$ and $z$ to length $2M$, as defined below.
\begin{equation*}
 \begin{aligned}
y_j &= x_{j}e^{-\pi ij^2\delta}   &\quad z_j &= e^{\pi i(j+s)^2\delta}      &\quad  \hspace{5mm}  \hbox{$0\leq j<M$}  \label{eq:l5}\\
y_j&=0   &\quad z_j&=e^{\pi i(j+s-2M)^2\delta}  &\quad  \hspace{5mm}
   \hbox{$M\leq j<2M$}
  \end{aligned}
\end{equation*}

\noindent
Taking into account the 2M-long sequence, the previous FRFT becomes
 \begin{equation}
 G_{k+s}(x,\delta)=e^{-\pi i (k+s)^2 \delta}\sum_{j=0}^{2M-1}\! y_{j}z_{k-j}=e^{-\pi i (k+s)^2 \delta}{DFT}_k^{-1}[{{DFT}(y){DFT}(z)}]    \label {eq:l6}
\end{equation}  
\noindent
Where $DFT$ is the Discrete Fourier Transform, and $DFT^{-1}$ is the inverse of $DFT$. For an M-long sequence $z$, we have 
 \begin{equation}
 DFT_{k}(z)=\sum_{j=0}^{M-1}z_{j}e^{-\frac{2\pi}{M}jk} \quad \quad  DFT^{-1}_{k}(z)=\frac{1}{M-1}\sum_{j=0}^{M-1}z_{j}e^{\frac{2\pi}{M}jk}  \label {eq:l7}
\end{equation}  This procedure is referred to in the literature as the Fast FRFT Algorithm with a total computational cost of $20Mlog_{2}M + 44M$ operations \cite {bailey1991fractional}.\\
\noindent
We assume that $\scrF[f](y)$ is zero outside the interval $[-\frac{a}{2}, \frac{a}{2}]$; $\beta=\frac{a}{M} $ is the step size of the $M$ input values of $\scrF[f](y)$, defined by $y_{j}=(j-\frac{M}{2}) \beta$ for $ 0 \leq j <M$. Similarly, $\gamma$ is the step size of the $M$ output values of $f(x)$, defined by $x_{k}=(k-\frac{M}{2}) \gamma$ for $ 0 \leq k <M$.\\
 By choosing the step size $\beta$ on the input side and the step size $\gamma$ in the output side, we fix the FRFT parameter $\delta=\frac{\beta\gamma}{2\pi}$ and yield \cite{nzokem2021fitting} the density function $f$ (\ref{eq:l7}) at $x_{k+s}$.
 \begin{equation*} 
 \begin{aligned}
   &f(x_{k+s}) = \frac{1}{2\pi}\int_{-\infty}^{+\infty}\! \scrF[f](y)e^{ix_{k+s}y} \mathrm{d}y \approx \frac{1}{2\pi}\int_{-a/2}^{a/2}\! \scrF[f](y)e^{ix_{k+s}y} \mathrm{d}y\\
   &\approx  \frac{\gamma}{2\pi}\sum_{j=0}^{M-1}\! \scrF[f](y_{j})e^{2\pi i(k+s-\frac{M}{2})(j-\frac{M}{2})\delta}=\frac{\gamma}{2\pi}e^{-\pi i(k+s-\frac{M}{2})M\delta}G_{k+s}(\scrF[f](y_{j})e^{-\pi i jM\delta},-\delta)
\end{aligned}
\end{equation*}
\noindent
 We have :
  \begin{align}
\tilde{f}(x_{k+s}) = \frac{\gamma} {2\pi}e^{-\pi i(k+s-\frac{M}{2})M\delta}G_{k+s}(\scrF[f](y_{j})e^{-\pi i jM\delta},-\delta) \hspace{5mm}
   \hbox{ $0\leq s <1$}
 \label {eq:l7}
 \end{align}  

\subsection{ FRFT of QN-long weighted sequence}
\noindent
The Q-point rule Composite Newton-Cotes Quadrature (\ref{eqp10}) is integrated into the Fractional Fourier (FRFT) algorithm (\ref{eq:l7}) to produce the FRFT of QN-long weighted sequence.\\
\noindent
We assume that $\scrF[f](x)$ is zero outside the interval $[-\frac{a}{2}, \frac{a}{2}]$, $M=QN$ and $\beta=\frac{a}{M} $ is the step size of the $M$ input values $\scrF[f](y)$, defined by $y_{j+Qp}=(Qp+ j - \frac{M}{2}) \beta$ for $0 \leq p <N$ and $0 \leq j <Q$. Similarly, the output values of $f(x)$ is defined by $x_{Ql+f+s}=(Ql+f+s-\frac{M}{2}) \gamma$ for $0 \leq l <N$, $0 \leq f <Q$ and $0\leq s\leq 1$.
 \begin{equation} 
 \begin{aligned}
f(x_{Ql+f+s})&=\frac{1}{2\pi}\int_{\infty}^{+\infty}e^{i y x_{Ql+f+s}} F[f](y)dy = \frac{1}{2\pi}\int_{-a/2}^{a/2}e^{i y x_{Ql+f+s}} F[f](y)dy\\
 &= \frac{1}{2\pi}\sum_{p=0}^{N-1} \int_{y_{Qp}}^{y_{Qp + Q}}e^{i y x_{Ql+f+s}} F[f](y)dy \hspace{2mm} \hbox{ (composite rule)} \label {eq:l8}
\end{aligned}
\end{equation}
\noindent
Based on the Lagrange interpolating integration over [$y_{Qp}$, $y_{Qp+Q}$]\cite{aubain2020, nzokem_2021}, we have:
\begin{equation} 
\int_{y_{Qp}}^{y_{Qp + Q}}{e^{i y x_{Ql+f+s}}F[f](y)}dy  \approx \beta\sum_{j=0}^{Q} w_{j} {e^{i y_{j + Qp} x_{Ql+f+s}} F[f](y_{j + Qp})}\label {eq:l9} 
 \end{equation}
\noindent
$\tilde{f}(x_{Ql+f+s})$ is  the approximation of $f(x_{Ql+f+s})$ and (\ref{eq:l8}) becomes
 \begin{equation} 
 \begin{aligned}
&\tilde{f}(x_{Ql+f+s}) =\frac{\beta}{2\pi}\sum_{p=0}^{N-1} \sum_{j=0}^{Q} w_{j} F[f](y_{j + Qp})e^{i x_{Ql+f+s} y_{j + Qp}} \hspace{2mm}
   \hbox{ $0\leq s <1$}\\
&= \frac{\beta}{2\pi}\sum_{j=0}^{Q} \sum_{p=0}^{N-1}{w_{j}F[f](y_{j + Qp})e^{2\pi i\delta (Ql+f+s-\frac{M}{2}) (Qp+ j - \frac{M}{2})}}  \hspace{2mm}
   \hbox{$\beta \gamma=2\pi\delta$}\\
&= \frac{\beta}{2\pi}e^{-\pi i\delta M(Ql+f+s-\frac{M}{2})} G_{Ql+f+s}(w_{j}\scrF[f](y_{j + Qp})e^{-\pi i(j + Qp)M\delta},-\delta)\label {eq:l10}
\end{aligned}
\end{equation}
\noindent
We have the following inverse Fourier transform function 
  \begin{align}
\resizebox{.9\hsize}{!}{$ \tilde{f}(x_{Ql+f+s})=\frac{\beta}{2\pi}e^{-\pi i\delta M(Ql+f+s-\frac{M}{2})} G_{Ql+f+s}(w_{j}\scrF[f](y_{j + Qp})e^{-\pi i(j + Qp)M\delta},-\delta) $}
 \label {eq:l11}
 \end{align} 
\noindent
FRFT is applied on the  QN-long weighted sequence $\left\{w_{j}\scrF[f](y_{j + Qp})e^{-\pi i(j + Qp)M\delta} \right \}_{\substack{0\leq j < Q \\ 0\leq p<N}}$.

\subsection{Composite of FRFTs: FRFT of Q-long weighted sequence FRFT of N-long sequence}
\noindent
We consider $\tilde{f}_{QN}(x_{Ql+f+s})$, the approximation of $f(x_{Ql+f+s})$ and the expression (\ref{eq:l10}) becomes
 \begin{equation*} 
 \begin{aligned}
&\tilde{f}_{QN}(x_{Ql+f+s}) =\frac{\beta}{2\pi}\sum_{p=0}^{N-1} \sum_{j=0}^{Q} w_{j} F[f](y_{j + Qp})e^{i x_{Ql+f+s} y_{j + Qp}}\\
&= \frac{\beta}{2\pi}\sum_{j=0}^{Q}w_{j} \sum_{p=0}^{N-1}{F[f](y_{j + Qp})e^{2\pi i\delta (Ql+f+s-\frac{M}{2}) (Qp+ j - \frac{M}{2})}} \hspace{5mm}
   \hbox{ $\beta \gamma=2\pi\delta$}\\
&= \frac{\beta}{2\pi}e^{-\pi i\delta M(Ql+f+s-\frac{M}{2})}\sum_{j=0}^{Q}w_{j}e^{2\pi i\delta (Ql+f+s-\frac{M}{2})j} \sum_{p=0}^{N-1}{F[f](y_{j + Qp})}e^{2\pi i\delta (Ql+f+s-\frac{M}{2})Qp}\\
&= \frac{\beta}{2\pi}e^{-\pi i\delta M(Ql+f+s-\frac{M}{2})}\sum_{j=0}^{Q}w_{j} G_{l+\frac{f+s}{Q}}(\xi_{p},-\delta Q^{2})e^{2\pi i\delta (Ql-\frac{M}{2})j}e^{2\pi i\delta (f+s)j}
\end{aligned}
\end{equation*}
\noindent
 $G_{l+\frac{f+s}{Q}}$ is a FRFT on N-long complex sequence $\{\xi_{p}\}_{0\leq p < N}$. Let us have $\alpha_{1}=-\delta Q^{2}$.
 \begin{equation} 
 \begin{aligned}
  \label {eq:l113}
G_{l+\frac{f+s}{Q}}(\xi_{p},\alpha_{1})&=\sum_{p=0}^{N-1}{\xi_{p}e^{-2\pi  i (l+\frac{f+s}{Q})p\alpha_{1}}} \quad \xi_{p}=e^{-\pi i M p Q \delta}F[f](y_{j + Qp})
\end{aligned}
\end{equation}
 $\tilde{f}_{QN}(x_{Ql+f+s})$ becomes
 \begin{equation} 
 \begin{aligned}
  \label {eq:l113a}
\tilde{f}_{QN}(x_{Ql+f+s}) = \frac{\beta}{2\pi}e^{-\pi i\delta M(Ql+f+s-\frac{M}{2})}G_{f+s}(z_{j}, -\delta)
\end{aligned}
\end{equation}
\noindent
$G_{f+s}$ is a FRFT on Q-long complex sequence  $\{z_{j}\}_{0\leq j \leq Q}$. Let us have  $\alpha_{2}=-\delta$.
 \begin{equation*} 
 \begin{aligned}
G_{f+s}(z_{j},\alpha_{2})&=\sum_{j=0}^{Q}z_{j}e^{-2\pi i(f+s)j\alpha_{2}} \quad \quad z_{j}=w_{j} G_{l+\frac{f+s}{Q}}(\xi_{p},\alpha_{1})e^{2\pi i\delta (Ql-\frac{M}{2})j}
\end{aligned}
\end{equation*}
\noindent
To sum up, we have the following expression for  $\tilde{f}_{QN}(x_{Ql+f+s})$
\begin{equation}
 \begin{aligned}
\label {eq:l115}
\resizebox{.9\hsize}{!}{$\tilde{f}_{QN}(x_{Ql+f+s})=\frac{\beta}{2\pi}e^{-\pi i\delta M(Ql+f+s-\frac{M}{2})}G_{f+s}(G_{l+\frac{f+s}{Q}}(\xi_{p},\alpha_{1})w_{j} e^{2\pi i\delta (Ql-\frac{M}{2})j},\alpha_{2})$}
\end{aligned}
\end{equation}

\begin{figure}[ht]
\vspace{-0.3cm}
    \centering
\hspace{-0.5cm}
  \begin{subfigure}[b]{0.32\linewidth}
    \includegraphics[width=\linewidth]{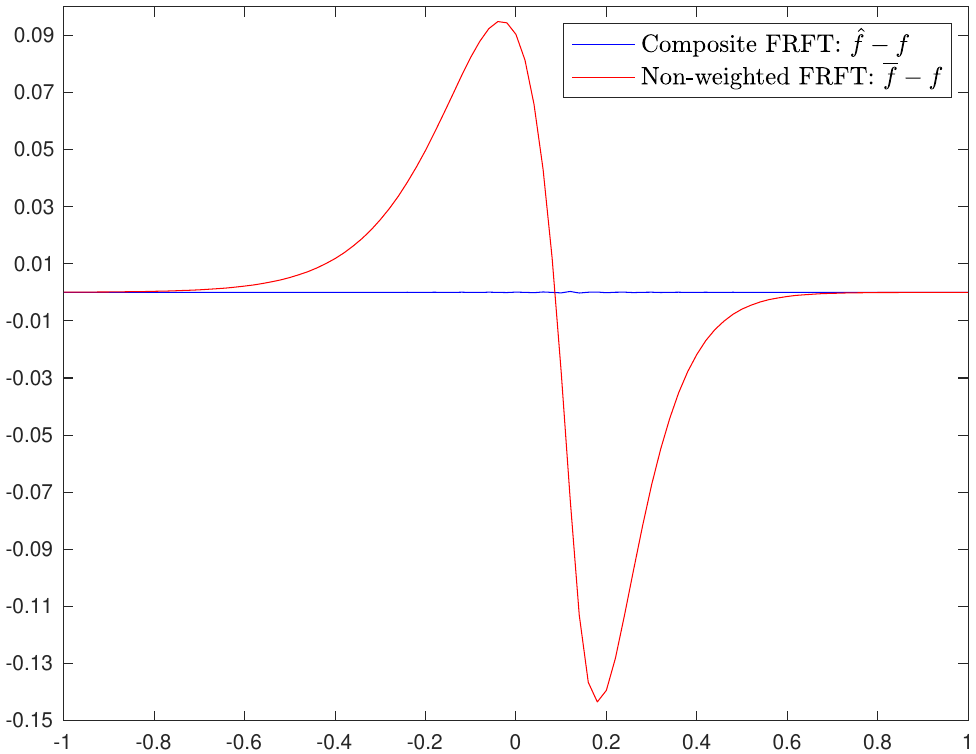}
\vspace{-0.6cm}
     \caption{$\textbf{Q=2}$}
         \label{fig21}
  \end{subfigure}
\hspace{-0.3cm}
  \begin{subfigure}[b]{0.32 \linewidth}
    \includegraphics[width=\linewidth]{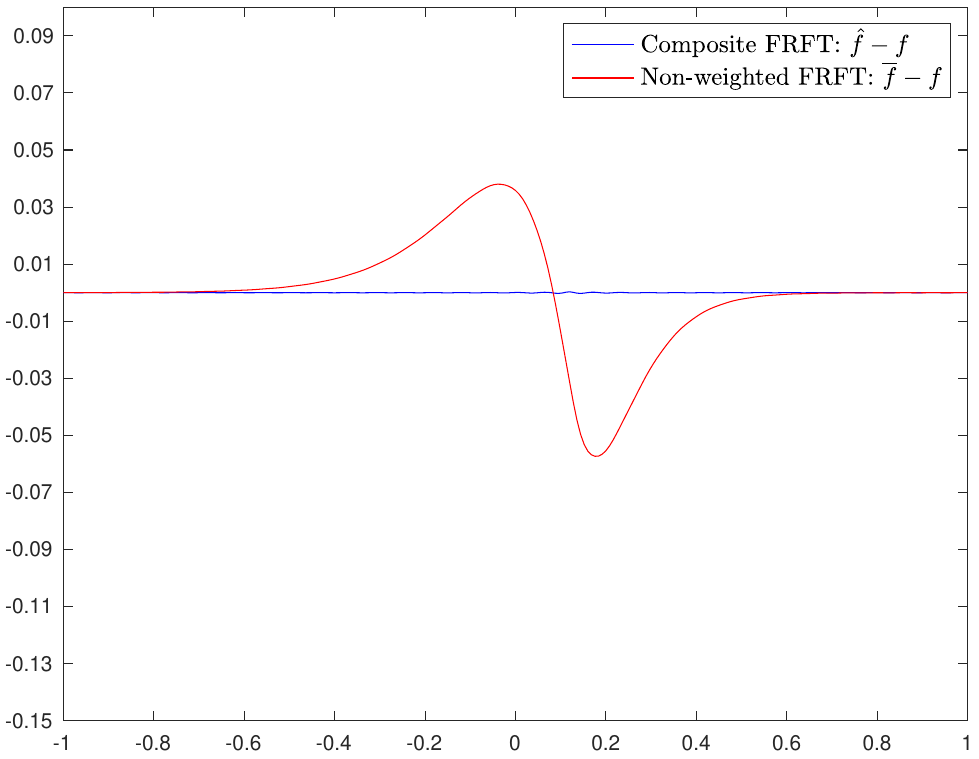}
\vspace{-0.6cm}
     \caption{$\textbf{Q=5}$}
         \label{fig22}
          \end{subfigure}
\hspace{-0.3cm}
  \begin{subfigure}[b]{0.32 \linewidth}
    \includegraphics[width=\linewidth]{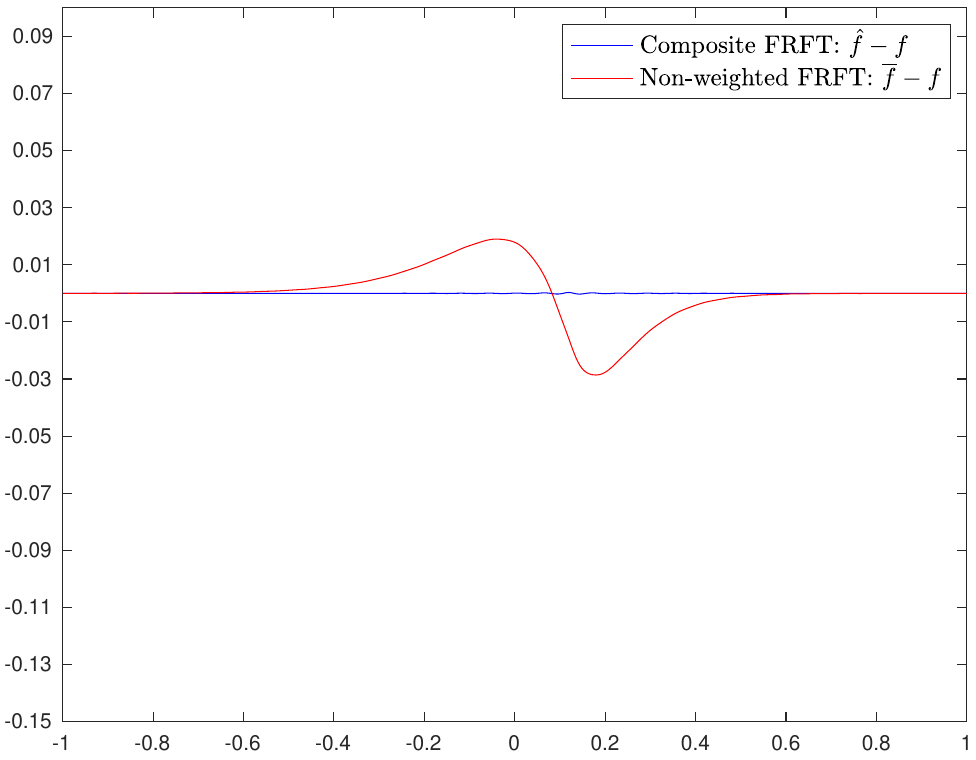}
\vspace{-0.6cm}
     \caption{$\textbf{Q=10}$}
         \label{fig23}
          \end{subfigure}
\vspace{-0.6cm}
  \caption{VG{*} Probability density error ($\textbf{$\hat{f}$ - $f$}$): $\textbf{N=5000}$ $\textbf{b-a=100}$}
  \label{fig2}
\vspace{-0.6cm}
\end{figure}

\noindent
By comparing the inverse Fourier transform the formula (\ref{eq:l11}) and (\ref{eq:l115}), we conclude that the FRFT ($G_{Ql+f}$) can be written as a composition of two FRFTs ($G_{f+s}$ and $G_{l+\frac{f+s}{Q}}$). 
\begin{equation}
 \begin{aligned}
 \label {eq:l116}
\resizebox{.9\hsize}{!}{$G_{Ql+f+s}(w_{j}\scrF[f](y_{j + Qp})e^{-\pi i(j + Qp)N\delta},\alpha_{2})=G_{f+s}\left(G_{l+\frac{f+s}{Q}}(\xi_{p},\alpha_{1})w_{j}e^{2\pi i\delta(Ql-\frac{M}{2})j},\alpha_{2} \right) $}
\end{aligned}
\end{equation}
\noindent
Using the Non-weighted FRFT(\ref{eq:l7}) and the composite FRFTs (\ref{eq:l116}), the risk-neutral probability density error of the Variance-Gamma $(\mu, \delta,\alpha,\theta,\sigma)$ model is estimated in Fig \ref{fig2}. The pattern is similar to the results in Fig \ref{fig1}.
\subsection{Composite of FRFTs: FRFT of N-long weighted sequence FRFT of Q-long sequence}
\begin{figure}[ht]
\vspace{-0.4cm}
    \centering
\hspace{-0.5cm}
  \begin{subfigure}[b]{0.32\linewidth}
    \includegraphics[width=\linewidth]{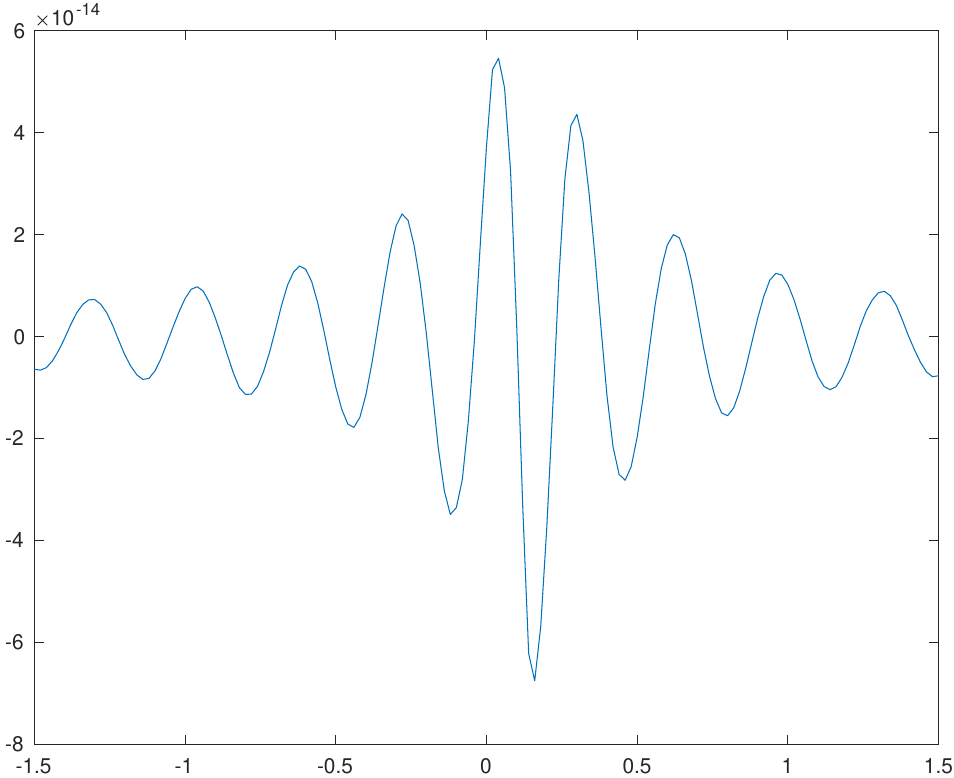}
\vspace{-0.6cm}
     \caption{$\textbf{Q=2}$}
         \label{fig21}
  \end{subfigure}
\hspace{-0.3cm}
  \begin{subfigure}[b]{0.32 \linewidth}
    \includegraphics[width=\linewidth]{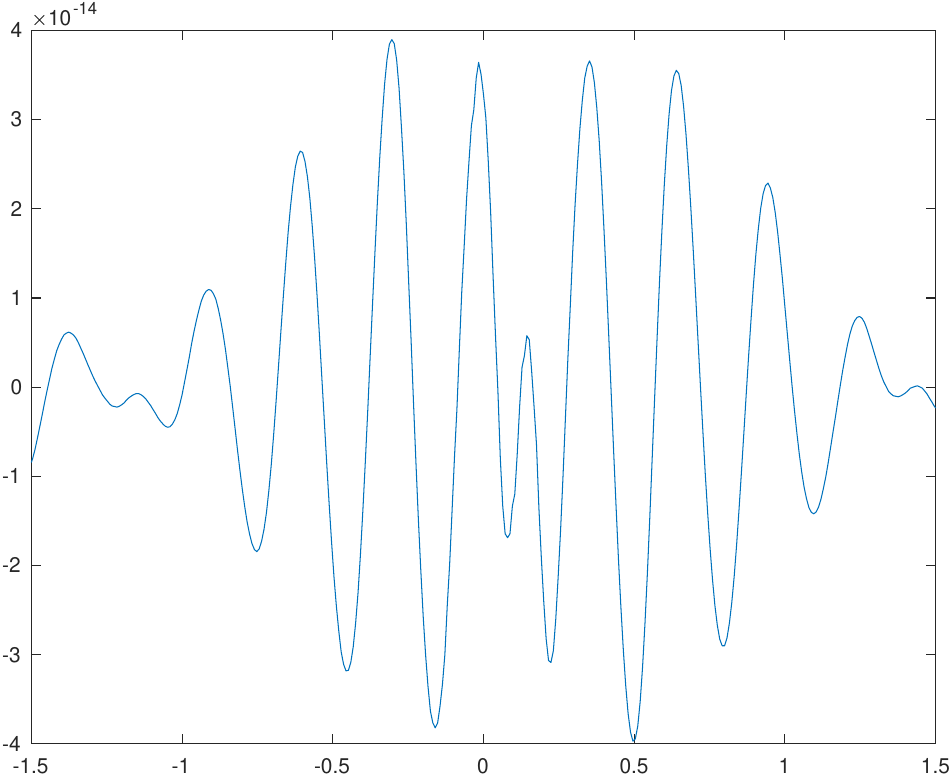}
\vspace{-0.6cm}
     \caption{$\textbf{Q=5}$}
         \label{fig22}
          \end{subfigure}
\hspace{-0.3cm}
  \begin{subfigure}[b]{0.32 \linewidth}
    \includegraphics[width=\linewidth]{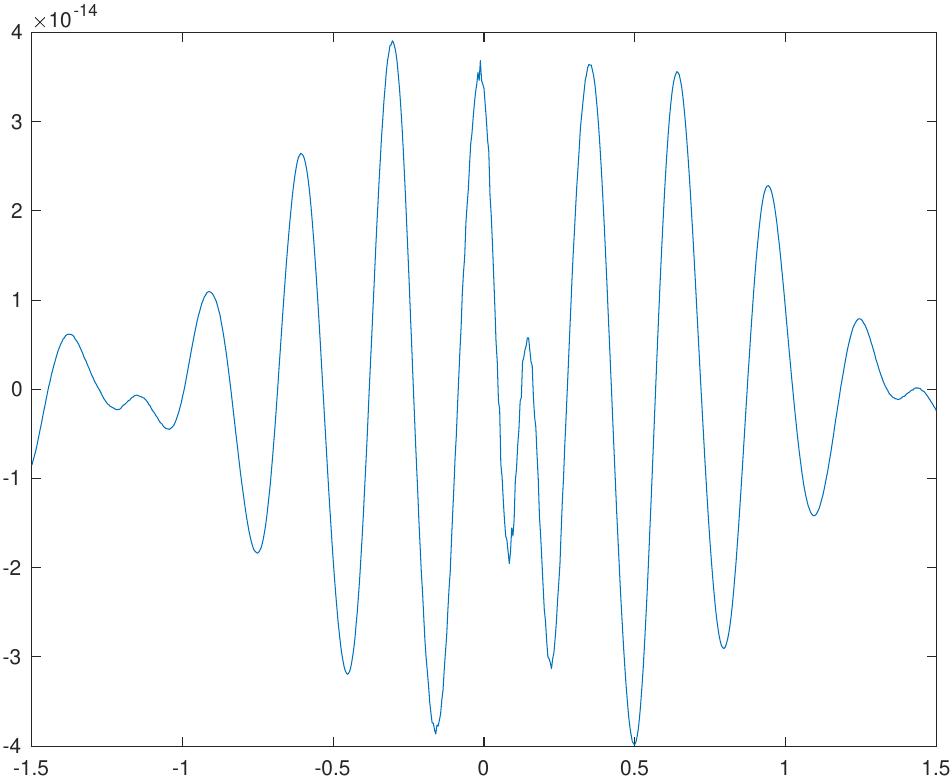}
\vspace{-0.6cm}
     \caption{$\textbf{Q=10}$}
         \label{fig23}
          \end{subfigure}
\vspace{-0.6cm}
  \caption{VG{*} Probability density error ($\textbf{$\tilde{f}_{NQ}(x_{k})$ - $\tilde{f}_{QN}(x_{k})$}$): $\textbf{N=5000}$ $\textbf{b-a=100}$}
  \label{fig25}
\vspace{-0.6cm}
\end{figure}

\noindent
The expression (\ref{eq:l10}) is recalled and we have
 \begin{equation*} 
 \begin{aligned}
&\tilde{f}_{NQ}(x_{Ql+f+s}) =\frac{\beta}{2\pi}\sum_{p=0}^{N-1} \sum_{j=0}^{Q} w_{j} F[f](y_{j + Qp})e^{i x_{Ql+f+s} y_{j + Qp}} \hspace{5mm}
   \hbox{ $\beta \gamma=2\pi\delta$} \\
&=\frac{\beta}{2\pi}e^{-\pi i\delta M(Ql+f+s-\frac{M}{2})}\sum_{p=0}^{N-1}\!e^{2\pi i\delta (Ql+f+s-\frac{M}{2})Qp}\sum_{j=0}^{Q}\!w_{j}{F[f](y_{j + Qp})e^{2\pi i\delta (Ql -\frac{M}{2})j}}e^{2\pi i\delta (f+s)j}\\
&=\frac{\beta}{2\pi}e^{-\pi i\delta M(Ql+f+s-\frac{M}{2})}\sum_{p=0}^{N-1}G_{f+s}(z_{j}, -\delta)e^{-\pi i\delta M Qp}e^{2\pi i(l+\frac{f+s}{Q})\delta Q^{2}p}
\end{aligned}
\end{equation*}
\noindent
$G_{f+s}$ is a fractional Fourier transform (FRFT) on Q-long complex sequence $\{z_{j}\}_{0\leq j < Q}$ 
\begin{equation} 
 \begin{aligned}
G_{f+s}(z_{j},\alpha_{2}) =\sum_{j=0}^{Q}z_{j}e^{-2\pi i(f+s)j\alpha_{2}} \quad \quad z_{j}=w_{j}F[f](y_{j + Qp})e^{2\pi i\delta (Ql -\frac{M}{2})j}
\end{aligned}
\end{equation}
$\tilde{f}_{NQ}((x_{Ql+f+s})$ becomes
 \begin{equation} 
 \begin{aligned}
\tilde{f}_{NQ}((x_{Ql+f+s}) = \frac{\beta}{2\pi}e^{-\pi i\delta M(Ql+f+s-\frac{M}{2})}G_{l+\frac{f+s}{Q}}(\xi_{p}, \alpha_{1})
\end{aligned}
\end{equation}
\noindent
$G_{l+\frac{f+s}{Q}}$ is a FRFT on the N-long complex sequence  $\{\xi_{p}\}_{0\leq p <{N}}$ 
 \begin{equation} 
 \begin{aligned}
G_{l+\frac{f+s}{Q}}(\xi_{p},\alpha_{1})=\sum_{p=0}^{N-1}{\xi_{p}e^{-2\pi  i (l+\frac{f+s}{Q})p\alpha_{1}}} \quad \quad \xi_{p}=G_{f+s}(z_{j},\alpha_{2})e^{-\pi i\delta M Qp}
\end{aligned}
\end{equation}
\noindent
As a result, we have the following expression for  $\tilde{f}_{NQ}(x_{Ql+f+s})$
\begin{equation}
 \begin{aligned}
 \label {eq:l120}
\resizebox{.9\hsize}{!}{$\tilde{f}_{NQ}(x_{Ql+f+s}) = \frac{\beta}{2\pi}e^{-\pi i\delta M(Ql+f+s-\frac{M}{2})}G_{l+\frac{f+s}{Q}}(G_{f+s}(z_{j}, -\delta)e^{-\pi i\delta M Qp}, -\delta Q^{2})$}
\end{aligned}
\end{equation}
\noindent
By comparing the inverse Fourier transform of the formula (\ref{eq:l11}) and (\ref{eq:l120}), we conclude that the FRFT ($G_{Ql+f+s}$)
is a composition of two FRFTs (($G_{l+\frac{f+s}{Q}}$) and ($G_{f+s}$ )). The result (\ref{eq:l121}) is similar to the result (\ref{eq:l116}).
\begin{equation}
 \begin{aligned}
 \label {eq:l121}
\resizebox{.9\hsize}{!}{$G_{Ql+f+s}(w_{j}\scrF[f](y_{j+Qp})e^{-\pi i(j + Qp)N\delta},\alpha_{2})=G_{l+\frac{f+s}{Q}}\left(G_{f+s}(z_{j},\alpha_{2})e^{-\pi i\delta M Qp},\alpha_{1}\right)$}
\end{aligned}
\end{equation}
\noindent
The numerical computation of (\ref{eq:l115}) and (\ref{eq:l120}) in Fig \ref{fig25} shows that the composite FRFTS in (\ref{eq:l116}) and (\ref{eq:l121}) are equal not only algebraically but also numerically.
\section{ILLUSTRATION EXAMPLES}
\subsection{Variance-Gamma VG Distribution}
\noindent
In the study, the VG model has five parameters: parameters of location ($\mu$), symmetric ($\delta$), volatility ($\sigma$), and the Gamma parameters of shape ($\alpha$) and scale ($\theta$). The VG model density function is proven to be (\ref {eq:l01}).
\begin{align}
f(y) =\frac {1} {\sigma\Gamma(\alpha) \theta^{\alpha}}\int_{0}^{+\infty}\!\frac{1}{\sqrt{2\pi \nu}}e^{-\frac{(y-\mu-\delta \nu)^2}{2\nu\sigma^2}}\nu^{\alpha -1}e^{-\frac{\nu}{\theta}} \,d\nu \label{eq:l01}
 \end{align}
\noindent
The VG model density function (\ref {eq:l01}) has an analytical expression with a modified Bessel function of the second kind. The expression can be obtained by making some transformations and changing variables in (\ref {eq:l01}).\\
\begin{align}
-\frac{(y-\mu-\delta \nu)^2}{2\nu\sigma^2}-\frac{\nu}{\theta} = \delta(\frac{y-\mu}{\sigma^2}) -\frac{1}{2\sigma^2}(\delta^2 + \frac{2\sigma^2}{\theta})\nu - (\frac{y-\mu}{\sigma})^2\frac{1}{2\nu}
 \end{align}

(\ref {eq:l01}) becomes 
\begin{align}
f(y) =\frac {e^{\delta(\frac{y-\mu}{\sigma^2})}} {\sqrt{2\pi}\sigma\Gamma(\alpha) \theta^{\alpha}}\int_{0}^{+\infty}\! e^{ -\frac{1}{2\sigma^2}(\delta^2 + \frac{2\sigma^2}{\theta})\nu - (\frac{y-\mu}{\sigma})^2\frac{1}{2\nu}}\nu^{\alpha -\frac{3}{2}} \,d\nu \label{eq:l03}
 \end{align}
\noindent
By considering the modified Bessel function of the second kind ($k_{\alpha}(z)$)\cite{NISTDLMF}, we have
\begin{align}
k_{\alpha}(z)=\frac{1}{2}(\frac{1}{2}z)^{\alpha}\int_{0}^{+\infty}\!e^{(-t -\frac{z^2}{4t})} \frac{1}{t^{\alpha +1}} \,dt \quad \quad  |arg(z)|\leq \frac{\pi}{4} \label{eq:l04}
 \end{align}
\noindent
$k_{\alpha}(z)$ is the second kind of solution for the modified Bessel's equation.
\begin{align}
z^2 \frac{d^{2}w}{dz^{2}} + z \frac{dw}{dz} - (z^2 + \alpha^2)w =0 \label{eq:l04a}
 \end{align}
\noindent
By changing variable, $u=\frac{1}{2\sigma^2}(\delta^2 + \frac{2\sigma^2}{\theta})\nu$, (\ref {eq:l03}) becomes 
\begin{align}
f(y)=\frac {2e^{\delta(\frac{y-\mu}{\sigma^2})}} {\sqrt{2\pi}\sigma\Gamma(\alpha) \theta^{\alpha}}\left(\frac{|y-\mu|}{\sqrt{\delta^2 + \frac{2\sigma^2}{\theta}}}\right)^{\alpha-\frac{1}{2}} k_{-\alpha +\frac{1}{2}}\left(\frac{\sqrt{\delta^2 + \frac{2\sigma^2}{\theta}}|y-\mu|}{\sigma^2}\right)  \label{eq:l05}
 \end{align}
\noindent
A special case $f(\mu)$ provide a simplified expression 
 \begin{equation} 
 \begin{aligned}
\resizebox{.9\hsize}{!}{$f(\mu) =\frac {1} {\sqrt{2\pi}\sigma\Gamma(\alpha) \theta^{\alpha}}\int_{0}^{+\infty}\!e^{ -\frac{1}{2\sigma^2}(\delta^2 + \frac{2\sigma^2}{\theta})\nu} \nu^{\alpha -\frac{3}{2}}\,d\nu =\frac{\Gamma(\alpha-\frac{1}{2})}{\sqrt{2\pi\theta}\sigma\Gamma(\alpha)( 1 + \frac{\theta}{2}\frac{\delta^2}{\sigma^2})^{\alpha -\frac{1}{2}}} \label{eq:l051}$}
 \end{aligned}
\end{equation}
\noindent
Table \ref {tab2} provides estimation results of the five parameters $(\mu,\delta,\alpha,\theta,\sigma)$ of the  Variance-Gamma variable. VG model estimation in the first raw of Table \ref {tab2} was obtained by Maximum likelihood method \cite{nzokem2021fitting, Nzokem_Montshiwa_2023, nzokem2024a}. The VG(*) model estimation in the second row of Table \ref {tab2} is the Equivalent Martingale Measure (EMM)  \cite{nzokem_2021, nzokem2023european, nzokem2022} and generates a risk-neutral probability density function for 0.8 year period. 
\begin{table}[ht]
\vspace{-0.3cm}
\caption{ VG Parameter Estimation}
\label{tab2} 
\vspace{-0.3cm}
\centering
\begin{tabular}{@{}lccccccc@{}}
\toprule
\textbf{Model} & \textbf{$\mu$} &\textbf{$\delta$} &  \textbf{$\sigma$} &  \textbf{$\alpha$} &   \textbf{$\theta$}\\  \midrule
\textbf{VG}&{$0.08476896$}&{$-0.0577418$}&{$1.02948292$}&{$0.88450029$}&{$0.93779517$} \\ 
\textbf{VG (*)}&{$0.11998901$}&{$-0.0343164$}&{$0.10294829$}&{$2.54736083$}&{$0.98780338$} \\ \bottomrule 
\end{tabular}
\vspace{-0.3cm}
\end{table} 

\begin{figure}[ht]
\vspace{-0.3cm}
    \centering
\hspace{-0.5cm}
  \begin{subfigure}[b]{0.32\linewidth}
    \includegraphics[width=\linewidth]{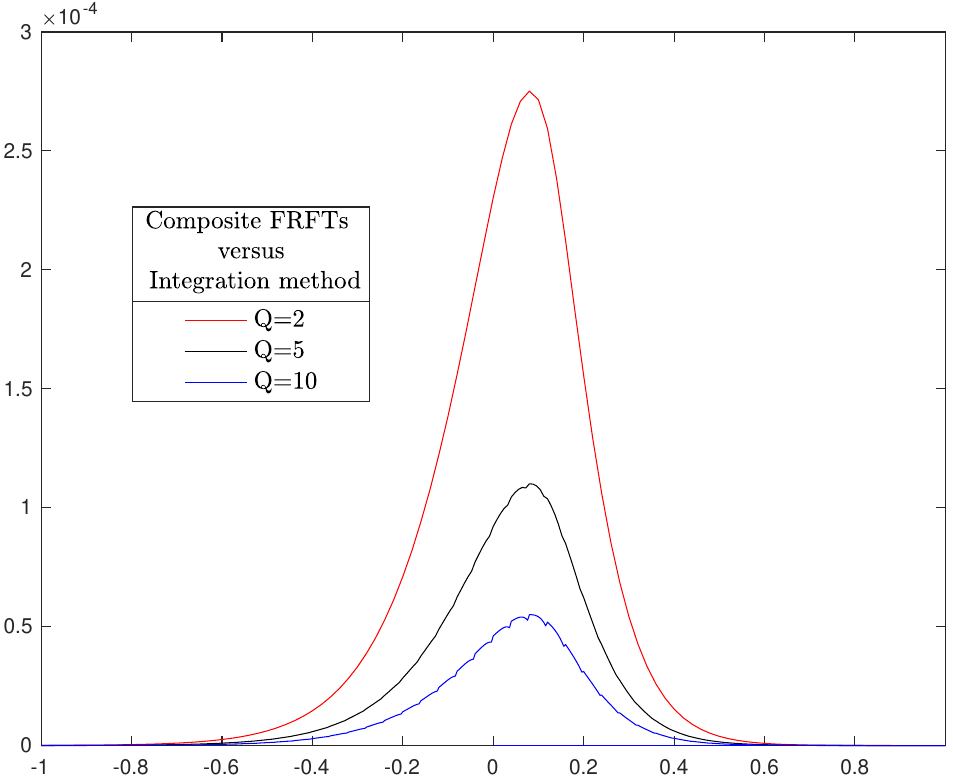}
\vspace{-0.3cm}
     \caption{VG{*} PDF Error ($\textbf{$\hat{f}$ - $\check{f}$}$) }
         \label{fig41}
  \end{subfigure}
\hspace{-0.3cm}
  \begin{subfigure}[b]{0.32\linewidth}
    \includegraphics[width=\linewidth]{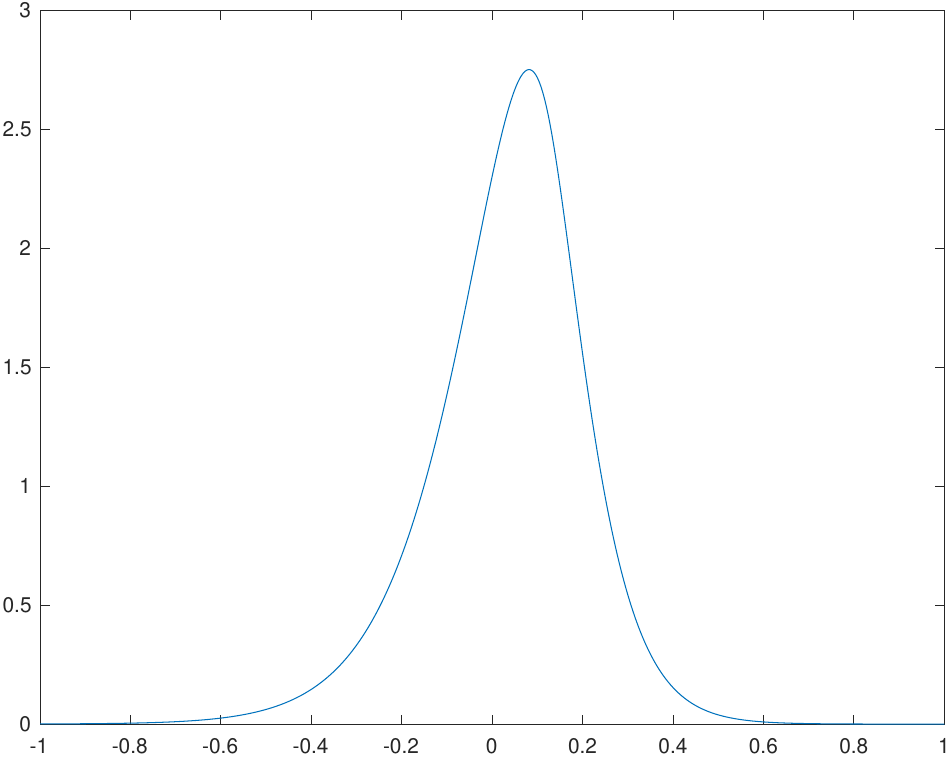}
\vspace{-0.3cm}
     \caption{VG{*} Density function}
         \label{fig42}
          \end{subfigure}
\hspace{-0.3cm}
  \begin{subfigure}[b]{0.32 \linewidth}
    \includegraphics[width=\linewidth]{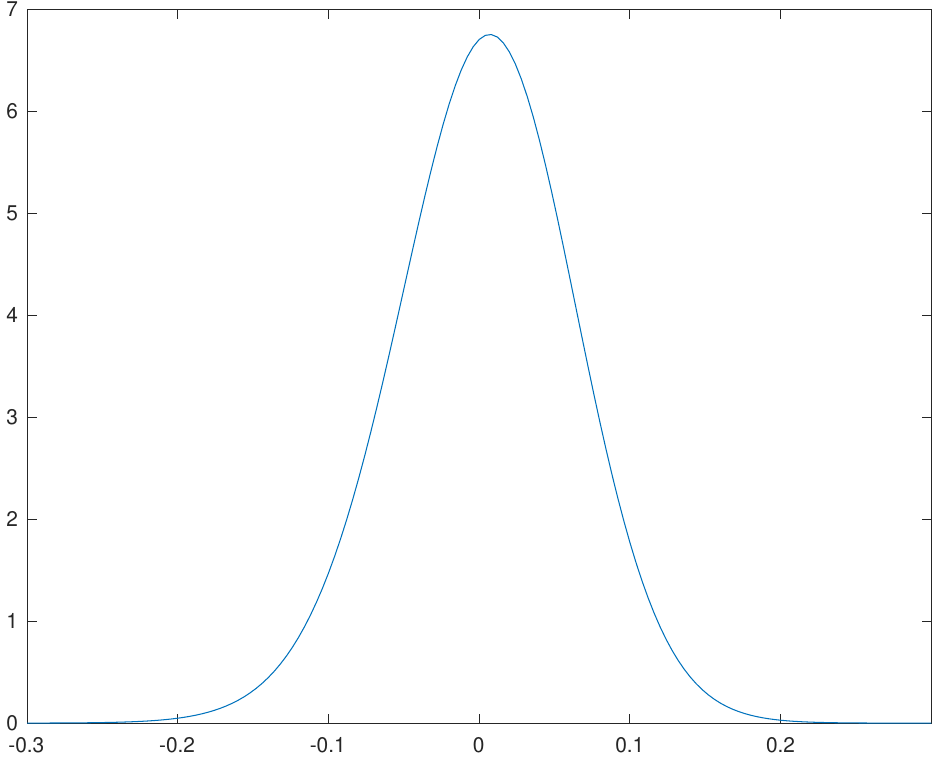}
\vspace{-0.3cm}
     \caption{GTS{*} Density function}
         \label{fig43}
          \end{subfigure}
\vspace{-0.3cm}
  \caption{ Risk Neutral Probability Density function (PDF) }
  \label{fig4}
\vspace{-1.0cm}
\end{figure}

\noindent
Fig \ref{fig42} provides a graphical representation of risk-neutral probability associated with the Variance Gamma (VG{*}) model for a 0.8-year period.\\
For special case: $f(\hat{\mu}) \approx 0.8552$ for VG model and  $f(\hat{\mu}) \approx 2.5949$ for VG{*} model.\\

\noindent
The Fourier transform function of the Variance Gamma distribution has an explicit closed form and the inverse function are summarised as follows. 
\begin{align}
 \scrF[f](x) &=\frac{e^{-i\mu x}}{\left(1+\frac{1}{2}\theta \sigma^{2}x^{2} + i\delta\theta x\right)^{\alpha}} \quad \quad f(y) = \frac{1}{2\pi}\int_{-\infty}^{+\infty}e^{iy x + \Psi(-x)}dx  \label {eq:l02}
  \end{align}

\noindent
Based on the integral method in proposition 1.1 (\ref{eqp10}), the Non weighted FRFT  and the composite FRFTs developed in (\ref{eq:l7}) and (\ref{eq:l115}) respectively,  we have the following numerical estimation methods of the inverse function (\ref{eq:l02}).
 \begin{align}
 \check{f}(x_{k}) &= \frac{\beta}{M} \sum_{p=0}^{N-1} \sum_{j=0}^{Q} W_{j} e^{(iy_{j+Qp} x_{k})}F[f](y_{j+Qp}) \hspace{5mm} \hbox{Integral approximation} \label {eq:l002}\\
 \overline{f}(x_{k}) &= \frac{\gamma} {2\pi}e^{-\pi i(k-\frac{N}{2})N\delta}G_{k}(\scrF[f](y_{j})e^{-\pi i jN\delta},-\delta)\hspace{5mm} \hbox{Non-weighted FRFT}  \label {eq:l003}\\
 \hat{f}(x_{Ql+f})&= \frac{\beta}{2\pi}e^{-\pi i\delta M(Ql+f-\frac{M}{2})}G_{f+s}(z_{j},-\alpha_{2}) \hspace{5mm} \hbox{Composite FRFTs} \label {eq:l004}
  \end{align}
\noindent
The VG{*} Probability density error ($\textbf{$\hat{f}$ - $\check{f}$}$) in Fig \ref{fig41} shows that the Composite FRFTs (\ref{eq:l004}) outperforms the  Newton-Cotes integration method (\ref{eq:l002}), but the error is not significantly different .
\subsection{Generalized Tempered Stable (GTS) Distribution}
\noindent
We consider a GTS variable $Y=\mu +X =\mu + X_{+} - X_{-} \sim GTS(\textbf{$\mu$}, \textbf{$\beta_{+}$}, \textbf{$\beta_{-}$}, \textbf{$\alpha_{+}$},\textbf{$\alpha_{-}$}, \textbf{$\lambda_{+}$}, \textbf{$\lambda_{-}$})$ with $X_{+}\sim TS(\textbf{$\beta_{+}$}, \textbf{$\alpha_{+}$},\textbf{$\lambda_{+}$}, \textbf{$\lambda_{-}$})$ and  $X_{-}\sim TS(\textbf{$\beta_{-}$}, \textbf{$\alpha_{-}$},\textbf{$\lambda_{-}$}, \textbf{$\lambda_{-}$})$. The characteristic exponent can be written \cite{nzokem2022fitting, nzokem2023european, mca29030044}
  \begin{align}
\resizebox{0.9\hsize}{!}{$\Psi(\xi)=\mu\xi i+\alpha_{+}\Gamma(-\beta_{+})\left((\lambda_{+}-i\xi)^{\beta_{+}}-{\lambda_{+}}^{\beta_{+}}\right)+\alpha_{-}\Gamma(-\beta_{-})\left((\lambda_{-}+\xi)^{\beta_{-}}-{\lambda_{-}}^{\beta_{-}}\right) $}\label {eq:l27}
  \end{align}
\noindent
Table \ref {tab3} presents the estimation results of the seven parameters $(\textbf{$\beta_{+}$}, \textbf{$\beta_{-}$}, \textbf{$\alpha_{+}$},\textbf{$\alpha_{-}$}, \textbf{$\lambda_{+}$}, \textbf{$\lambda_{-}$})$ Generalized Tempered Stable (GTS) Distribution. GTS model estimation in the first raw of Table \ref {tab2} was obtained by Maximum likelihood method \cite{nzokem2022fitting, Nzokem_Montshiwa_2023, jrfm17120531}. 
\begin{table}[ht]
\vspace{-0.3cm}
\caption{ GTS Parameter Estimation}
\label{tab3} 
\vspace{-0.3cm}
\centering
\begin{tabular}{@{}lccccccc@{}}
\toprule
\textbf{Model} & \textbf{$\mu$} & \textbf{$\beta_{+}$} & \textbf{$\beta_{-}$} & \textbf{$\alpha_{+}$} & \textbf{$\alpha_{-}$}  & \textbf{$\lambda_{+}$}  & \textbf{$\lambda_{-}$}  \\ \midrule
\textbf{GTS} & -0.693477 & 0.682290 & 0.242579 & 0.458582 & 0.414443 & 0.822222 & 0.727607  \\ 
\textbf{GTS{*}} & -0.208043 & 0.682290 & 0.242579 & 0.594234 & 4.068436 & 84.667097 & 70.31591 \\ \bottomrule
\end{tabular}
\vspace{-0.6cm}
\end{table} 

\noindent
The GTS(*) model estimations in the second row of Table \ref {tab2} is the Equivalent Martingale Measure (EMM)  \cite{nzokem2023european} and generate a risk-neutral probability density function for 0.8 year period. Fig \ref{fig43} provides a graphical representation of risk-neutral probability associated with the GTS Distribution for a 0.8-year period.\\

\begin{figure}[ht]
\vspace{-0.4cm}
    \centering
\hspace{-0.5cm}
  \begin{subfigure}[b]{0.32\linewidth}
    \includegraphics[width=\linewidth]{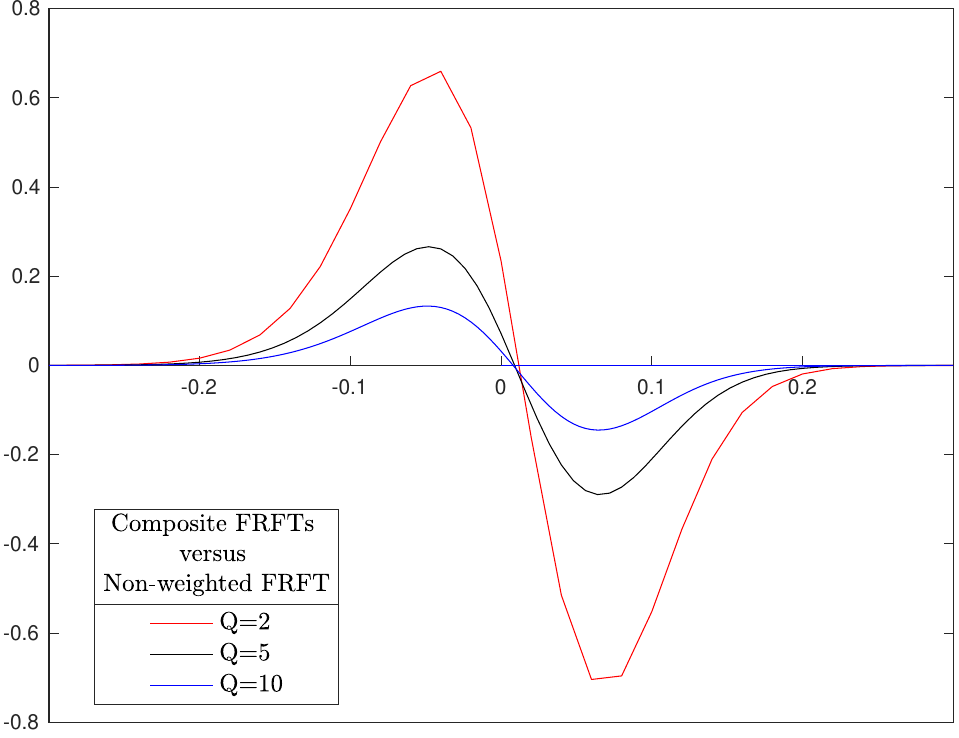}
\vspace{-0.6cm}
     \caption{$\textbf{$\check{f}(x_{k}) - \overline{f}(x_{k})$}$}
         \label{fig51}
  \end{subfigure}
\hspace{-0.3cm}
  \begin{subfigure}[b]{0.32\linewidth}
    \includegraphics[width=\linewidth]{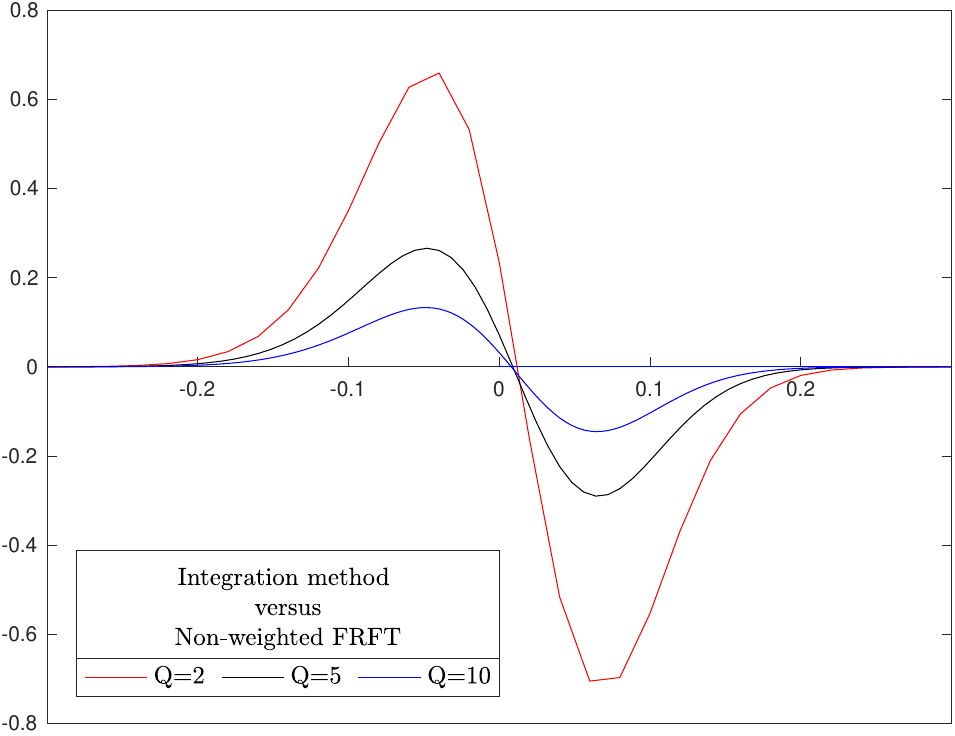}
\vspace{-0.6cm}
     \caption{ $\textbf{$\hat{f}(x_{k}) - \overline{f}(x_{k})$}$}
         \label{fig52}
          \end{subfigure}
\hspace{-0.3cm}
  \begin{subfigure}[b]{0.3\linewidth}
    \includegraphics[width=\linewidth]{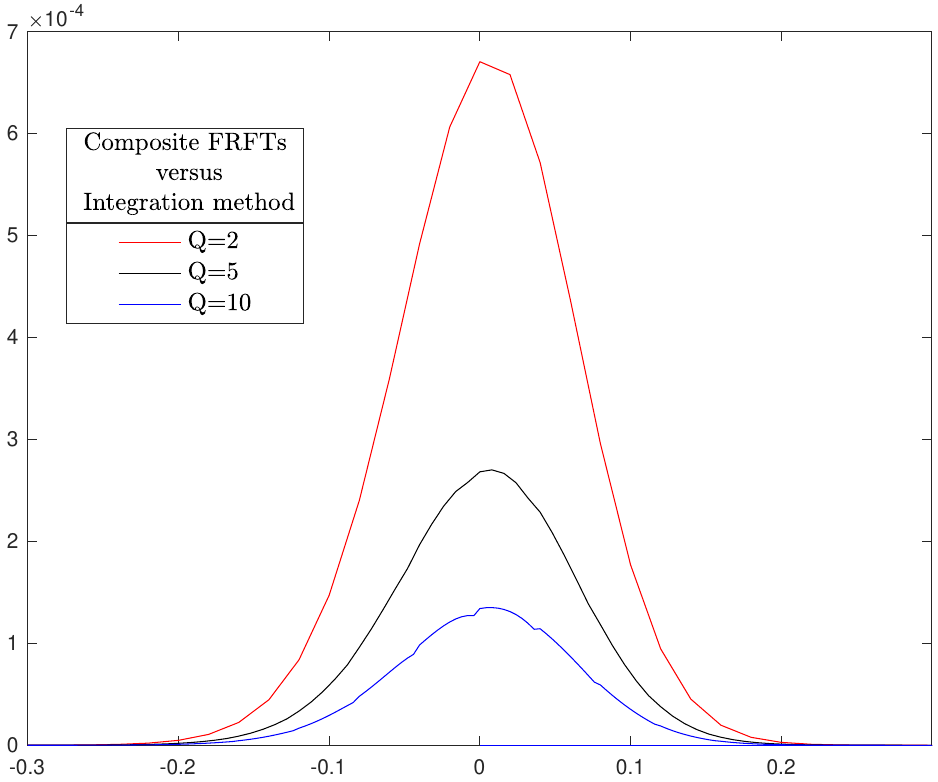}
\vspace{-0.6cm}
     \caption{ $\textbf{$\hat{f}(x_{k}) -  \check{f}(x_{k})$}$}
         \label{fig52}
          \end{subfigure}
\vspace{-0.3cm}
  \caption{GTS{*} Probability density error}
  \label{fig5}
\vspace{-0.8cm}
\end{figure}

\noindent
The characteristic function of the GTS variable, the Fourier Transform ($F(f)$), and the density function ($f$) have the following expression
\begin{equation}
 \vartheta(\xi)=E\left[e^{i Y \xi}\right]=e^{\Psi(\xi)} \quad \quad F[f](\xi) = \vartheta(-\xi) \quad f(y) = \frac{1}{2\pi}\int_{-\infty}^{+\infty}e^{iy x + \Psi(-x)}dx \label {eq:l2}
\end{equation}
 \noindent
The GTS probability density function has neither closed form nor analytic expression. Fig \ref{fig5}  provide an estimation of the GTS{*} risk neutral  Probability density error using the integration method ($\check{f}$) in (\ref{eq:l002}), the Non weighted FRFT ($\overline{f}$) in (\ref{eq:l002})and  the composite FRFTs ($\hat{f}$) in ( \ref{eq:l004}). The composite FRFTs and the integration method outperform the Non-weight FRFT in Fig \ref{fig51} and \ref{fig52}.

\section {Conclusion} 
The closed composite Newton-Cotes quadrature rule and the FRFT algorithm are reviewed in the paper. Both schemes are combined to yield a Composite of FRFTs. It is shown that the FRFT of a QN-long weighted sequence can be written as a composite of two FRFTs: the FRFT of a Q-long weighted sequence and the FRFT of an N-long sequence, and the FRFT of an N-long weighted sequence and the FRFT of a Q-long sequence. It is shown that the composite FRFTs have the commutative property  and work both algebraically and numerically. The composite of  FRFTs scheme was applied to analyse the probability density function error of the Variance-Gamma VG $(\mu, \delta,\alpha,\theta,\sigma)$ distribution and the Generalized Tempered Stable (GTS)(\textbf{$\beta_{+}$}, \textbf{$\beta_{-}$}, \textbf{$\alpha_{+}$},\textbf{$\alpha_{-}$}, \textbf{$\lambda_{+}$}, \textbf{$\lambda_{-}$}) Distribution. The results suggest that the composite of FRFTs outperforms both the simple non-weighted FRFT and the Newton-Cotes integration method, but the difference is less significant for the integration method.

\section*{Acknowledgement}
Notably, 2025 is dedicated to the 70$^{th}$ anniversary of academician Evgeny Evgenievich Tyrtyshnikov, whose contributions continue to inspire progress in numerical mathematics.

\bibliography{nzokem_pricing}

\end{document}